\input amstex
\documentstyle{amsppt}
\magnification 1200
\vcorrection{-9mm}
\input epsf

\topmatter
\title
        Quasipositive links and connected sums
\endtitle
\author
        S.~Yu.~Orevkov
\endauthor
\address
Steklov Mathematical Institute of Russian Academy of Sciences, Moscow.
\endaddress
\abstract
We prove that the connected sum of two links is quasipositive if and only if
each summand is quasipositive. The proof is based on the filling disk technique.
\endabstract
\thanks
This work was supported by the Russian Science Foundation under grant 19-11-00316
\endthanks
\endtopmatter

\def\refBK      {1} 
\def\refBe      {2} 
\def\refBKL     {3}
\def\refBO      {4} 
\def\refDe      {5} \let\refDemailly=\refDe
\def\refEli     {6} 
\def\refEr      {7} 
\def\refEV      {8} % Etnyre - Van Horn-Morris
\def\refHay     {9}
\def\refK       {10}
\def\refOrMM    {11}
\def\refOrBirk  {12}
\def\refRuTop   {13}
\def\refRuIII   {14}
\def\refRuSlice {15}
\def\refTa      {16} \let\refTakeuchi=\refTa

\def\eqX {1}

\def\thConnSum        {1}
\def\thSplitLink      {2}
\def\thSplitBraid     {3}
\def\conjConnSum      {4} \let\conjConnSumBraid=\conjConnSum
\def\remQPMM          {5}
\def\remConnSum       {6}
\def\remSplitSum      {7}
\def\remCombinat      {8}
\def\remSplitBraid    {9}

\def\sectQ   {2}
\def\lemQ    {10}

\def\sectSQ  {3}
\def\lemSQ   {11}

\def\C{\Bbb{C}}
\def\CP{\Bbb{CP}}
\def\eps{\varepsilon}
\def\sh{\operatorname{sh}}
\def\lk{\operatorname{lk}}

\document

\head 1. Introduction
\endhead

An $n$-braid is called {\it quasipositive} if it is a product of conjugates
of the standard (Artin's) generators $\sigma_1,\dots,\sigma_{n-1}$ of the
braid group $B_n$. A braid is called {\it strongly quasipositive} if it
is a product of braids of $\sigma_{j,k+1}=\tau_{k,j}\sigma_j\tau_{k,j}^{-1}$
for $j\le k$ where $\tau_{k,j}=\sigma_{k}\sigma_{k-1}\dots\sigma_j$.
Such braids are called {\it band generators} (they are also known as the
generators in the Birman-Ko-Lee presentation of $B_n$, see [\refBKL]).

All links in this paper are assumed to be oriented links in the 3-sphere $S^3$.
A link is called {\it (strongly) quasipositive} if it is the braid closure of a
(strongly) quasipositive braid. This terminology was introduced by Lee Rudolph
(see [\refRuTop, \refRuIII]) and now it has become standard in the knot theory.

The main result of this note is the (a) statement of the following theorem.

\proclaim{ Theorem \thConnSum } Let $L=L_1\# L_2$ be the connected sum of two links in $S^3$.
Then:

\smallskip
(a). $L$ is quasipositive if and only if $L_1$ and $L_2$ are. % quasipositive.

\smallskip
(b). $L$ is strongly quasipositive if and only if $L_1$ and $L_2$ are.
% strongly quasipositive.
\endproclaim

Note that the (b) statement of this theorem is an almost immediate
corollary of the main result of [\refRuIII] (see \S\sectSQ).
I added it just for the sake of completeness as well as
the following two theorems.

\proclaim{ Theorem \thSplitLink } 
Let $L=L_1\sqcup L_2$ be the split sum of two links in $S^3$.
Then:

\smallskip
(a). $L$ is quasipositive if and only if $L_1$ and $L_2$ are.

\smallskip
(b). $L$ is strongly quasipositive if and only if $L_1$ and $L_2$ are.
\endproclaim

Let $\sh_m:B_n\to B_{m+n}$ be the homomorphism defined on the generators by
$\sigma_k\mapsto \sigma_{m+k}$ (the $m$-shift).
If links $L_1$ and $L_2$ are represented by braids
$X_1\in B_m$ and $X_2\in B_n$, then $L_1\sqcup L_2$ and $L_1\#L_2$
can be represented by the braids
$$
    X_1 \sh_m(X_2)\in B_{m+n}\qquad\text{and}\qquad
    X_1 \sh_{m-1}(X_2)\in B_{m+n-1}     \eqno(\eqX)
$$
respectively. So, the next result is
a braid-theoretic counterpart of Theorem \thSplitLink.

\proclaim{ Theorem \thSplitBraid }
Let $X_1\in B_m$ and $X_2\in B_n$. Let $X=X_1 \sh_m(X_2)\in B_{m+n}$. Then:

\smallskip
(a). {\rm([\refOrBirk, Thm.~3.2])}
$X$ is quasipositive if and only if $X_1$ and $X_2$ are.

\smallskip
(b). $X$ is strongly quasipositive if and only if $X_1$ and $X_2$ are.
\endproclaim

\proclaim{ Conjecture \conjConnSum }
Let $X_1\in B_m$ and $X_2\in B_n$. Let $X=X_1 \sh_{m-1}(X_2)\in B_{m+n-1}$.
Then $X$ is quasipositive if and only if $X_1$ and $X_2$ are.
% quasipositive.
\endproclaim

\medskip\noindent
{\bf Remark \remQPMM.}
A particular case of Conjecture \conjConnSumBraid\ is the main result of
[\refOrMM] which states
that an $n$-braid $X$ is quasipositive if and only if the $(n+1)$-braid
$X\sigma_n$ is.

 \medskip\noindent
 {\bf Remark \remConnSum.}
 In [\refOrMM, Question~1] I asked if the minimal braid index representative
 of a quasipositive link is necessarily a quasipositive braid. In virtue of
 Theorem \thConnSum(a), an affirmative answer to this question 
 combined with arguments from [\refHay] would imply
 Conjecture \conjConnSum.
 Indeed, the self-linking number (the algebraic length minus the number of strings)
 of a quasipositive braid is maximal over 
 all braids representing the same link type, see e.g.~[\refRuSlice].
 It follows that, in the setting of Conjecture~\conjConnSum,
it is maximal for each $X_j$, $j=1,2$. Then, as shown in the proof
 of [\refHay, Thm.~1.2], $X_j$ can be transformed into a braid $X_j'$ with the minimal
 number of strings by conjugations and positive (de)stabilizations only.
 Hence Theorem \thConnSum(a) combined with an affirmative answer to
 [\refOrMM, Question~1] would imply the
 quasipositivity of $X_j'$ which is equivalent to that of $X_j$ by [\refOrMM, Thm.~1].
 See also [\refEV, \refHay] for some interesting results related to
 [\refOrMM, Question~1].

\medskip\noindent
{\bf Remark \remSplitSum.}
In particular, Theorem \thSplitBraid(a) implies that
{\it an $n$-braid is quasipositive, if so is the
$m$-braid given by the same braid word for some $m\ge n$}
(see [\refOrBirk, Thm.~3.1]). This fact is also an immediate
formal consequence from the invariance of quasipositivity under destabilizations
(see [\refOrMM] and Remark \remQPMM).
Indeed, if $X\in B_n$ and the $(n+1)$-braid given by the same
braid word is quasipositive, then evidently so is $X\sigma_n\in B_{n+1}$ whence,
by [\refOrMM], $X$ as well.

\medskip\noindent
{\bf Remark \remCombinat.}
All the discussed statements concerning (not strongly) quasipositive braids
are purely combinatorial whereas their proofs are based on the (almost) complex
analysis, PDE, etc. The only particular case where I know a combinatorial proof
is the statement in Remark \remSplitSum; see [\refOrBirk,~\S3.3].

\medskip\noindent
{\bf Remark \remSplitBraid.}
Conjecture \conjConnSum\ is a braid-theoretical counterpart of
Theorem~\thConnSum(a). Using the Birman-Ko-Lee version of Garside's theory [\refBKL],
one can easily prove several analogs of Theorem~\thConnSum(b) for braids.
However, they do not seem to be of any interest because the strong quasipositivity
is not invariant under conjugations. 

\medskip
\subhead Acknowledgements \endsubhead
I am grateful to Michel Boileau for attracting my attention to this topic
and to Nikolay Kruzhilin and Stefan Nemirovski for useful discussions.
Also I thank the referee for correcting some errors.

%%%%%%%%%%%%%%%%%%%%%%%%%%%%%%%%%%%%%%%%%%%%%%%%%%%%%%%%%%%%%%%%%%%
%%%%%%%%%%%%%%%%%%%%%%%%%%%%%%%%%%%%%%%%%%%%%%%%%%%%%%%%%%%%%%%%%%%
%%%%%%%%%%%%%%%%%%%%%%%%%%%%%%%%%%%%%%%%%%%%%%%%%%%%%%%%%%%%%%%%%%%
%%%%%%%%%%%%%%%%%%%%%%%%%%%%%%%%%%%%%%%%%%%%%%%%%%%%%%%%%%%%%%%%%%%
%%%%%%%%%%%%%%%%%%%%%%%%%%%%%%%%%%%%%%%%%%%%%%%%%%%%%%%%%%%%%%%%%%%
%%%%%%%%%%%%%%%%%%%%%%%%%%%%%%%%%%%%%%%%%%%%%%%%%%%%%%%%%%%%%%%%%%%
%%%%%%%%%%%%%%%%%%%%%%%%%%%%%%%%%%%%%%%%%%%%%%%%%%%%%%%%%%%%%%%%%%%

\head\sectQ. Quasipositive links: proof of the {\rm(a)} cases of the theorems
\endhead

The proofs of the (a) cases of the theorems are inspired by
Eroshkin's paper [\refEr] which is based on the filling disk technique
[\refBK, \refEli, \refK]. In fact, Theorems \thSplitLink(a)
and \thSplitBraid(a) are almost immediate consequences from [\refEr]
(see [\refOrBirk, Thm.~3.2] for more details).
Our proof of Theorem \thConnSum(a) is a combination of
a more precise version of Bedford-Klingenberg-Kruzhilin's Theorem
(with the smoothness of the Levi-flat hypersurface),
the smoothing of pseudoconvex domains, and the result from [\refBO]
which states that the boundary link of an algebraic curve in a
pseudoconvex ball is a quasipositive link.

\proclaim{ Lemma \lemQ }
Let $\Omega$ be a bounded pseudoconvex domain in $\C^2$ with smooth
boundary. Let $B\subset\Omega$ be a smooth Levi-flat hypersurface transverse
to $\partial\Omega$, and $A$ be a smooth two-dimensional surface transverse to
both $B$ and $\partial\Omega$. Suppose that $\Omega\setminus B$ has two
connected components $\Omega_1$ and $\Omega_2$.
Then for each $j=1,2$, there exists a strictly pseudoconvex domain
$\Omega_j^-\subset\Omega_j$ with smooth boundary such that the pairs
$(\Omega_j,A\cap\Omega_j)$ and $(\Omega_j^-,A\cap\Omega_j^-)$
are homeomorphic.
\endproclaim

\demo{ Proof } We identify $\C^2$ with an affine chart of $\CP^2$. Let $j=1$ or $2$.
For $z\in\Omega_j$, we define $d(z)$ as the distance from $p$ to $\partial\Omega_j$
with respect to the Fubini-Study metric on $\CP^2$. By [\refTa, Thm.~1],
$h:=-\log d$ is a plurisubharmonic function on $\Omega_j$.
Let $\varphi$ be a smooth function in $\C^2$ supported
on the unit ball and positive on it,
and let $\varphi_\eps(z)=\varphi(z/\eps)/\eps^4$
(thus $\varphi_\eps$ tends to a delta-function as $\eps\to0$).
Let $U_\eps=\{p\in\Omega_j\mid d(p)>\eps\}$.
Let $h_\eps$ be the convolution $h*\varphi_\eps$, it is smooth and
plurisubharmonic on $U_\eps$ (see e.g. [\refDe, Thm.~5.5]).
Then, for $a\gg1$ and $0<\eps\ll\exp(-a)$, the domain
$\Omega_j^-=\{z\in\C^2\mid h_\eps(z)+\eps\|z\|^2 < a\}$
has the required properties.
%[\refDe, Thm. (5.5)-(5.6)] (smoothing the corners)
\qed
\enddemo

\demo{ Proof of Theorem \thConnSum(a) } ($\Leftarrow$) Follows from (\eqX).

($\Rightarrow$) Suppose that $L$ is a quasipositive link in $S^3$ which
we identify with the unit sphere in $\C^2$.
By [\refRuTop] we may assume that $L=S^3\cap A$ where $A$ is an algebraic curve
transverse to $S^3$.
Let $\Gamma\subset S^3$ be a smooth embedded 2-sphere which defines the
decomposition of $L$ into the connected sum $L_1\# L_2$. This means that
$\Gamma$ divides $S^3$ into two 3-balls $B_1$ and $B_2$, and there is an
embedded segment $I\subset\Gamma$ such that $(L\cap B_j)\cup I=L_j$, $j=1,2$.
It follows from [\refBK, Thm.~1] that, after perturbing $\Gamma$ if necessary,
one can find a smoothly embedded Levi-flat $3$-ball $B$ transverse
\footnote{ The transversality is not stated in [\refBK, Thm.~1], however, 
           the proof of the smoothness of $B$ is nothing else
           than a proof of its transversality to $S^3$.}
to $S^3$ and such that $\partial B=\Gamma$.
Moreover (see also [\refEli]), there exists a smooth function $F:B\to\Bbb R$
with non-vanishing gradient, whose restriction to $\Gamma$ (we denote it by $f$)
is a Morse function, and all whose level surfaces are unions of holomorphic disks
and (for some critical levels) isolated points.

Let $\{p_1,p_2\}=\partial I=\Gamma\cap L$, we number the points $p_1$
and $p_2$ so that $f(p_1)\le f(p_2)$. Then the linking numbers of $L$ with
the level lines of $f$ are:
$$
  \lk\big(L,f^{-1}(c)\big) = \cases 1 & \text{ if $f(p_1)\le c\le f(p_2)$,} \\
                                    0 & \text{ otherwise.} \endcases
$$
Hence $F^{-1}(c)$ meets $A$ transversally at a single point when
$f(p_1)\le c\le f(p_2)$, and $F^{-1}(c)\cap A=\varnothing$ otherwise.
Therefore $B\cap A$ is an unknotted arc in $B$ with the
endpoints $p_1,p_2\in\Gamma=\partial B$
(it is crucial here that the fibers $F^{-1}(c)$ are unions of
disks because otherwise an arc cutting each fiber at most once might be knotted).

Let $\Omega_1$ and $\Omega_2$ be the two domains
into which $B$ divides the unit ball in $\C^2$.
Let $j=1$ or $2$. 
We have $\partial\Omega_j = B\cup B_j$, and the arc $I$ is isotopic to
$B\cap A$ in $B$ relative to the boundary whence the
homeomorphisms (see Figure~1):
$$
  (S^3,L_j) = (S^3,(B_j\cap A)\cup I) \cong
  (\partial\Omega_j, (B_j\cap A)\cup I) \cong
  (\partial\Omega_j, \partial(\Omega_j\cap A)).
$$
By Lemma \lemQ, we may approximate $\Omega_j$ by a strictly
pseudoconvex domain $\Omega_j^-$ with smooth boundary diffeomorphic to the $3$-sphere.
Then, by Eliashberg's result [\refEli, Thm.~5.1],  $\Omega_j^-$ is diffeomorphic to
the $4$-ball,
and then the quasipositivity
of $(\partial\Omega_j, \partial(\Omega_j\cap A))$ follows from [\refBO, Thm.~2].
\qed\enddemo

\vbox{
\epsfysize=30mm\epsfbox{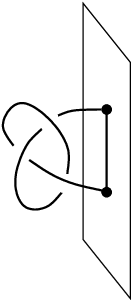}
\hskip -30pt\lower -10pt\hbox{$B_j$ \hskip 13pt $B_{3-j}$}
\hbox to 7.5em{
  \vbox{\pagewidth{5em}
    replace\par\noindent\hskip2pt
    $B_{3-j}$ by $B$\par\noindent
    $-\!\!-\!\!-\!\!-\!\!-\!\!-\!\!-\!\!{\to}$
    \vskip 5ex
  }
}
\epsfysize=30mm\epsfbox{a.eps}
\hskip -30pt\lower -10pt\hbox{$B_j$ \hskip 13pt $B$}
\hbox to 8em{
  \vbox{\pagewidth{5.5em}\noindent
    $I\sim A\cap B$\par\noindent\hskip 0.5em inside $B$
    \par\noindent
    $-\!\!-\!\!-\!\!-\!\!-\!\!-\!\!{\to}$
    \vskip 5ex
  }
}
\epsfysize=30mm\epsfbox{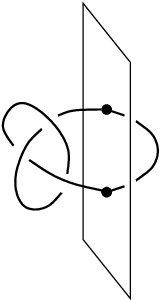}
\hskip -37pt\lower -10pt\hbox{$B_j$ \hskip 13pt $B$}
\medskip
\centerline{\bf Figure~1.}
}% end vbox
\medskip

%%%%%%%%%%%%%%%%%%%%%%%%%%%%%%%%%%%%%%%%%%%%%%%%%%%%%%%%%%%%%%%%%%%
%%%%%%%%%%%%%%%%%%%%%%%%%%%%%%%%%%%%%%%%%%%%%%%%%%%%%%%%%%%%%%%%%%%
%%%%%%%%%%%%%%%%%%%%%%%%%%%%%%%%%%%%%%%%%%%%%%%%%%%%%%%%%%%%%%%%%%%
%%%%%%%%%%%%%%%%%%%%%%%%%%%%%%%%%%%%%%%%%%%%%%%%%%%%%%%%%%%%%%%%%%%
%%%%%%%%%%%%%%%%%%%%%%%%%%%%%%%%%%%%%%%%%%%%%%%%%%%%%%%%%%%%%%%%%%%
%%%%%%%%%%%%%%%%%%%%%%%%%%%%%%%%%%%%%%%%%%%%%%%%%%%%%%%%%%%%%%%%%%%
%%%%%%%%%%%%%%%%%%%%%%%%%%%%%%%%%%%%%%%%%%%%%%%%%%%%%%%%%%%%%%%%%%%

\head\sectSQ. Strongly quasipositive links: proof of the {\rm(b)} cases of the theorems
\endhead

\proclaim{ Lemma \lemSQ }
Let $L$ be $L_1\# L_2$ or $L_1\sqcup L_2$.
Let $S$ be a Seifert surface of $L$ of maximal Euler characteristic.
Then the sphere $S^3$ can be presented as a union of embedded 3-balls
$S^3=B_1\cup B_2$ such that:
\roster
\item"(i)" $B_1\cap B_2 = \partial B_1=\partial B_2$;
\item"(ii)" $\partial(S\cap B_j)$ is $L_j$ for $j=1,2$;
\item"(iii)" $S\cap\partial B_1$ is an embedded segment if $L=L_1\# L_2$
and empty if $L=L_1\sqcup L_2$.
\endroster
\endproclaim

\demo{ Proof }
The lemma follows from the standard arguments used in the proof of the
additivity of the knot genus. Namely,
let $S^3=B_1\cup B_2$ be a splitting of $S^3$ involved in the definition
of the split or connected sum. We assume that $\partial B_1$ is transverse
to $S$. If $(B_1,B_2)$ is not as required, we choose a closed curve $C$ in
$S\cap \partial B_1$ which bounds a disk $D$ in $\partial B_1\setminus S$.
If $C$ bounds a disk $D'$ in $S$, then
$D\cup D'$ bounds a 3-ball
$B\subset B_j$, $j=1$ or $2$, and we may replace
$B_{3-j}$ with a thickening of $B_{3-j}\cup B$ which reduces the number
of components of $S\cap\partial B_1$.
Otherwise we attach a 2-handle to $S$ along $D$
and remove a closed component if it appears. This operation increases $\chi(S)$
which contradicts its maximality.
\qed\enddemo

If an $n$-braid $X$ is a product of $c$
band generators, then its braid closure admits a Seifert surface of
a special form which has Euler characteristic $n-c$.
It is obtained by attaching $c$ positively half-twisted
bands to $n$ parallel disks so that a band corresponding to $\sigma_{i,j}$
connects the $i$-th disk with the $j$-th disk; see details in [\refRuIII].
Such a surface is called a {\it quasipositive Seifert surface}.

\demo{ Proof of Theorems \thConnSum(b) and \thSplitLink(b)}
($\Leftarrow$) Follows from (\eqX).

$(\Rightarrow)$
Suppose that $L$ is strongly quasipositive. Let $S$ be a quasipositive
Seifert surface for $L$. By Bennequin's inequality [\refBe, Ch. II, Thm. 3]
the Euler characteristic of $S$ is maximal among all Seifert surfaces of $L$.
Hence the sphere $S^3$ can be cut into two 3-balls $S^3=B_1\cup B_2$ as
in Lemma \lemSQ.
Then $B_j\cap S$ is a Seifert surface of $L_j$ which is a full
subsurface of $S$ (see the definition in [\refRuIII, p.~231]),
thus $B_j\cap S$ is isotopic to a quasipositive surface by the main
result of [\refRuIII], whence the strong quasipositivity of $L_j$. \qed
\enddemo

\demo{ Proof of Theorem \thSplitBraid(b) } ($\Leftarrow$) Evident.

($\Rightarrow$) Let $S$ be the quasipositive Seifert surface. By the same
reasons as in the proof of Theorems \thConnSum(b) and \thSplitLink(b),
$S$ is a disjoint union $S_1\cup S_2$ with $\partial S_j=L_j$. Then each $S_j$ is a
quasipositive Seifert surface by construction. \qed
\enddemo

%%%%%%%%%%%%%%%%%%%%%%%%%%%%%%%%%%%%%%%%%%%%%%%%%%%%%%%%%%%%%%%%%%%
%%%%%%%%%%%%%%%%%%%%%%%%%%%%%%%%%%%%%%%%%%%%%%%%%%%%%%%%%%%%%%%%%%%
%%%%%%%%%%%%%%%%%%%%%%%%%%%%%%%%%%%%%%%%%%%%%%%%%%%%%%%%%%%%%%%%%%%
%%%%%%%%%%%%%%%%%%%%%%%%%%%%%%%%%%%%%%%%%%%%%%%%%%%%%%%%%%%%%%%%%%%
%%%%%%%%%%%%%%%%%%%%%%%%%%%%%%%%%%%%%%%%%%%%%%%%%%%%%%%%%%%%%%%%%%%
%%%%%%%%%%%%%%%%%%%%%%%%%%%%%%%%%%%%%%%%%%%%%%%%%%%%%%%%%%%%%%%%%%%
%%%%%%%%%%%%%%%%%%%%%%%%%%%%%%%%%%%%%%%%%%%%%%%%%%%%%%%%%%%%%%%%%%%

\enddocument